\documentclass{amsart}
\usepackage{amscd} 
\usepackage{amsfonts} 
\usepackage{amssymb} 
\usepackage{latexsym} 
\input{diagrams}
\newcommand{\ncm}{\newcommand}

\def\M{\mathcal{M}}

\newtheorem{theorem}{Theorem}[section]
\newtheorem{prop}[theorem]{Proposition}
\newtheorem{lemma}[theorem]{Lemma}
\newtheorem{cor}[theorem]{Corollary}
\newtheorem{lem&def}[theorem]{Lemma \& Definition}
\newtheorem{definition}[theorem]{Definition}

\ncm{\End}{\mbox{\rm End}\,}

\def\Hom{\mbox{\rm Hom}\,}

\def\id{\mbox{\rm id}}

\def\into{\hookrightarrow}
\def\to{\rightarrow}
\def\Ind{\mbox{\rm Ind}}
\def\Res{\mbox{\rm Res}}
\def\Coind{\mbox{\rm CoInd}}

\def\o{\otimes}    %tensor product 
\ncm{\rarr}[1]{\stackrel{#1}{\longrightarrow}}
\ncm{\larr}[1]{\stackrel{#1}{\longleftarrow}}

\def\eps{\varepsilon}

\def\du1{\hat 1}

\def\-1{_{(-1)}}
\def\0{_{(0)}}
\def\1{_{(1)}}
\def\2{_{(2)}}
\def\3{_{(3)}}

\def\|{\, | \,}

\def\du1{\hat 1}

\def\ract{\triangleleft}

\begin{document}

\title{Centralizers and Inverses to Induction as Equivalence of Categories}
\author{Lars Kadison}
\date{\today}
\address{Matematiska Institutionen \\ G{\" o}teborg 
University \\ 
S-412 96 G{\" o}teborg, Sweden} 
\email{lkadison@c2i.net} 
\thanks{16D90, 18E05, 20D25, 22D30}
\subjclass{}

\begin{abstract} 
Given a ring homomorphism $B \rightarrow A$, consider its centralizer $R = A^B$, 
 bimodule endomorphism ring $S = \End {}_BA_B$ and sub-tensor-square ring $T = (A \o_B A)^B$. 
Nonassociative tensoring by the cyclic modules $R_T$ or ${}_SR$  leads to an  equivalence of categories inverse 
to the functors of induction of restricted $A$-modules or  restricted coinduction of $B$-modules
in case $A \| B$  is separable, H-separable, split or left depth two (D2).  If $R_T$ or ${}_SR$ are
 projective, this property characterizes separability or splitness for a ring extension.  Only in
the case of H-separability is $R_T$ a progenerator, which replaces the key module $A_{A^e}$ for an
Azumaya algebra $A$.  After establishing these characterizations, we characterize left D2 extensions in terms of the module $T_R$, and ask whether a weak generator condition
on $R_T$ might characterize left D2 extensions as well, possibly a problem in $\sigma(M)$-categories or
its generalizations.
  We also show that the centralizer of a depth two extension is a normal subring in the sense of Rieffel
and pre-braided commutative.  
For example, its normality yields a Hopf subalgebra analogue of a factoid for subgroups and their centralizers,
and a special case of a conjecture that D2 Hopf subalgebras are normal.  
\end{abstract} 
\maketitle

\section{Introduction and Preliminaries}

Given a ring homomorphism $B \rightarrow A$, we pass to its induced bimodule ${}_BA_B$ and
define its centralizer $R = A^B$, 
 bimodule endomorphism ring $S = \End {}_BA_B$ and sub-tensor-square ring $T = (A \o_B A)^B$. 
The ring $T$ (defined in eq.~(\ref{eq: tee}) below) replaces the enveloping ring $A^e$ which plays a major role
in the study of separability and Azumaya algebras \cite{DI}.  We will see that it plays a similar role
in the study of separability and H-separability below, perhaps in a troubling way from the point of view
of Morita equivalence and its generalizations. When we  
tensor in a nonassociative way by the cyclic module $R_T$,  an induced $A$-module $A \o_B M$ of a separable extension becomes  naturally $A$-isomorphic to   $M$ itself (Lemma~\ref{lem-sep}).  This shows that induction of $A$-modules
is an equivalence of categories which characterizes separability if
$R_T$ is just projective (Theorem~\ref{th-sep}) and characterizes H-separability if $R_T$
is a generator (Theorem~\ref{th-hsep}).  Thus only in the case of H-separability is this
 equivalence part of a Morita equivalence, in this case tensoring by the progenerator
module $R_T$ with endomorphism ring isomorphic to the center of $A$.  

 A split extension is the dual of a separable extension in a unusual way involving the endomorphism ring \cite{K}; if we now pass from $T$ to its $R$-dual $S$ in depth two theory \cite{KS},
we might experiment with the cyclic module ${}_SR$ (which does not generally form a bimodule with the right $T$-module action). It turns out that tensoring the restricted coinduction of a $B$-module $N$ by ${}_SR$ is naturally isomorphic
to $N$ (Lemma~\ref{lem-split}). Thus restricted coinduction
is an equivalence from the category of modules $\M_B$ to its image subcategory,
a  characterization of split extension if ${}_SR$ is projective
(Theorem~\ref{th-split}).   
 
Left depth two extensions have a Galois theory and duality theorems for
actions and coactions explored in \cite{KS,LK, LK2, LK3, CQR}.  What we do instead in this paper
is show that a left D2 extension $A \| B$ has properties similar to a split extension as well as a separable extension.
In Theorems~\ref{th-ld2} and~\ref{th-left d2 coind} we show that a necessary condition for $A \| B$ to be left D2 is that 
restricted coinducted restriction and inducted restriction are equivalences as functors from the category
of $A$-modules with inverse from the image subcategory defined by $- \o_S R$ and $R \o_T -$ 
once again.  Given the known examples of D2 and separable extensions
independently scattered on a Venn diagram  on the one hand, and the characterization
of H-separability (a strong form of both
D2 extension and separable extension) on the other, we ask whether a weak generator condition on $R_T$ together with the equivalence data might characterize left D2 extensions, perhaps with a flatness condition in addition.    
This might have an answer in the  the theory of $\sigma(M)$-categories or its generalizations by Wisbauer \cite{W} and Dress
a subject we bring to bear in section~5 on a related question.

    In section~4 of this paper we revisit a notion of normal subring in the sense of Rieffel \cite{R}. We work with the simplest possible definition of normality in \cite{R} since this already
coincides with the notion of normality for Hopf subalgebras (Proposition~\ref{prop-hop}). It is also
a good notion for stating that the centralizer of a depth two extension is a normal subring (Proposition~\ref{prop-nc}, 
and the closely related fact that it is pre-braided commutative in section~5).  For example, this yields a Hopf subalgebra analogue of a factoid for subgroups and their centralizers (Proposition~\ref{prop-crop}),
and a special case of a conjecture that D2 Hopf subalgebras are normal (Corollary~\ref{cor-tys}).

In a final section we discuss further categorical considerations and results such as which ring extensions
generalizing separable extensions as well as D2 extensions possess the property that the inverse of induction-restriction of modules is tensoring by the centralizer module $R_T$.
We also discuss a depth two condition for bimodules using bicategories and endomorphism rings, and
a functorial characterization of left D2 extension $A \| B$ in terms of induction from $B$ to $A$
on the one hand and induction from $R$ to $T$, or coinduction from $R$ to $S$, on the other hand.

\subsection{Preliminaries for the general reader} 
The starting point in this paper is a ring homomorphism $g: B \to A$ of associative unital rings preserving unity
and admitting the possibility of noncommutative rings. (We may equally well work with $K$-algebras over a commutative
ring $K$, $K$-symmetric bimodules and $K$-linear morphisms throughout. In this paper then unadorned tensors between rings are over
the integers.)  
 The homomorphism induces a natural bimodule
${}_BA_B$ via $b \cdot a \cdot b' = g(b)ag(b')$ if $g: B \to A$, as well as bimodules ${}_AA_B$
and ${}_BA_A$ defined via $g$ on one side only. The same holds for an $A$-module ${}_AM$;
it is a $B$-module ${}_BM$ via $g$. This defines the usual functor $R$ (of ``restriction'') from the category of
left $A$-module denoted by ${}_A\M$ with $A$-linear morphism into
the category ${}_B\M$.  
The data we focus on is contained
in the subring $g(B) \subseteq A$, leading us to suppress $g$ and 
 view $B \to A$ as a  \textit{ring extension} 
$A \| B$, a \textit{proper ring extension} if $B \into A$ is monic.   Sometimes a condition on $A_B$, such as being a generator module, 
entails that $A \| B$ is a proper extension.  
 If $B$ is commutative and $g$ maps
$B$ into the center $Z(A)$ of $A$, then $A \| B$ is the $B$-algebra $A$, where a proper ring extension
is a faithful algebra. 

  If ${}_AM_A$ denotes an $A$-$A$-bimodule, we let $$ M^B := \{ m \in M | bm = mb, \, \forall \, b \in B \}. $$
Note that this is isomorphic to the group of $A$-$A$-bimodule homomorphisms from the tensor-square to $M$, 
\begin{equation}
\label{eq: em}
 \Hom ({}_A A \o_B A_A, {}_AM_A) \stackrel{\cong}{\longrightarrow} M^B  \ \ \  F \longmapsto F(1_A \o 1_A) 
\end{equation}
with inverse $m \longmapsto (x \o y \mapsto xmy)$ for $m \in M^B, x,y \in A$.  

In particular, we let $M = A$ and note the centralizer subring of $A$, 
\begin{equation}
\label{eq: are}
R := A^B \cong \Hom ({}_AA \o_B A_A, {}_AA_A),
\end{equation}
via $r \mapsto (x \o y \mapsto xry)$.

As another particular case, we let $M = A \o_B A$, the \textit{tensor-square} with its usual endpoint $A$-$A$-bimodule structure, so that the endomorphism ring
(under composition)
\begin{equation}
\End {}_AA \o_B A_A \cong (A \o_B A)^B := T
\end{equation}
induces the following ring structure on the construction $T$ for a ring extension $A \| B$, where 
 we suppress a possible summation over simple tensors and use a Sweedler-like
notation to write $t = t^1 \o t^2 \in T \subseteq A \o_B A$.
\begin{equation}
\label{eq: tee}
tu = u^1 t^1 \o t^2 u^2, \ \ \ 1_T = 1_A \o 1_A
\end{equation}
We will fix the notation $T$ throughout the paper for this ring construction while
not making explicit its dependence on the extension $A \| B$. Note the left ideal $(A \o_B A)^A$ of 
\textit{Casimir tensor elements}.  

From the $\End$-representation of $T$ we note that $A \o_B A$ is the left $T$-module:
\begin{equation}
{}_T (A \o_B A): \ \ \ t \cdot (x \o y) = xt^1 \o t^2 y 
\end{equation}
for all $t \in T, x,y \in A$.  

Thus any module ${}_AM$ (restricted then) induced to $A \o_B M$ becomes a left $T$-module
via the canonical isomorphism $A \o_B M \cong A \o_B A \o_A M$:
\begin{equation}
{}_T(A \o_B M): \ \ \ t \cdot (a \o m) = a t^1 \o t^2 m 
\end{equation}
for $m \in M, a \in A, t \in T$. 

 Given the bimodule ${}_AM_A$, the $\End$-representation of $T$ 
and eq.~(\ref{eq: em}) leads via composition from the right to
the module $(M^B)_T$ given by
\begin{equation}
\label{eq: embeetee}
(M^B)_T: \ \ \ m \cdot t = t^1 m t^2 \ \ (m \in M^B, t \in T).
\end{equation}

In particular, we obtain the module $R_T$, important in this paper, defined by
\begin{equation}
R_T: \ \ \ r \cdot t = t^1 r t^2 \ \ \ (r \in R, t \in T).
\end{equation}
This module was studied in \cite{KK} as a generalized Miyashita-Ulbrich module
with roots in Hopf algebra and group representations.  

The bimodule ${}_RT_R$ is important in the depth two theory in e.g.\ \cite{KS,LK,LK2,LK3}
where $T$ acts as a right bialgebroid over $R$ ,
 and plays a minor role in this paper:
\begin{equation}
{}_RT_R: \ \ \ r \cdot t \cdot r' := rt^1 \o t^2 r' 
\end{equation}
for $t \in T, r,r' \in R$, a restriction of scalars to a submodule of the tensor-square $A$-$A$-bimodule. 

We also need the construction $S := \End {}_BA_B$, dual to $T$ as a left $R$-bialgebroid in the
depth two theory \cite{KS,LK,LK2}. In the ring $S$ under the usual composition we note the
right ideal $\Hom ({}_BA_B, {}_BB_B)$. For any $B$-module $N_B$,
we note that the \textit{coinduced} right $A$-module $\Hom (A_B, N_B)$
(given by $(ha)(x) = h(ax)$) is also a right $S \o B$-module; i.e., with commuting
right $S$- and $B$-action since $(hb)(\alpha(a)) =  (h \circ \alpha)(ba)$ for
$\alpha \in S$, $b \in B$.  In other words, the $A$-module action restricted to a $B$-module action commutes with
the $S$-module action to give a module $\Hom ({}_BA, {}_BN)_{S \o B}$
(tensor over the integers). 

Since each $r \in R$ defines an endomorphism in $S$ in two ways - by left multiplication or by right multiplication,
\begin{equation}
\lambda_r(x) = rx, \ \ \ \rho_r(x) = xr \ \ \ (x \in A, r \in R)
\end{equation}
which commute $\lambda_r \rho_{s} = \rho_{s} \lambda_r$, we work  with an $R$-$R$-bimodule structure on $S$
given by
\begin{equation}
{}_RS_R: \ \ \ r \cdot \alpha \cdot s =   \lambda_r \rho_s \alpha = r \alpha(-) s
\end{equation}
which is the appropriate bimodule for the left bialgebroid theory in \cite{KS,LK, LK2, LK3}.  

Finally, there is the natural left $S$ action on $A$ by evaluation, which restricted to
the centralizer becomes
\begin{equation}
{}_SR: \ \ \ \alpha \cdot r = \alpha(r)
\end{equation}
which is also in $R$  (for $\alpha \in S, r\in R$). The module 
${}_SR$, studied to an extent in \cite{LK}, will be important to our study of split extensions in the next section
and left D2 extensions in the following section.

%%%%%%%%%%%%%%%%%%%%%%%%%%%%%%%%%%%%%%%%%%%%%%%%%%%%%%%%%%%%%%%%%%%%%%%%%%%%%%%%%%%%%%%%%%%%%%%%%%%%%%%%%%%%%%%%%%%%%
\section{Separable and Split Extensions}

Recall that a ring extension $A \| B$ is \textit{separable} if the mapping
$\mu: A \o_B A \to A$ given by multiplying components, $\mu(x \o y) = xy$, is a split $A$-$A$-epimorphism;
an element $e \in \mu^{-1}(1_A)$ is called a separability element and characterizes
separable extensions with its two properties, $ae^1 \o e^2 = e^1 \o e^2a$ for all $a \in A$ and $e^1 e^2 = \mu(e) = 1$.  

\begin{lemma}
\label{lem-sep}
Given  a separable extension $A \| B$ and a module ${}_AM$,
the mapping $r \o_T (a \o_B m) \stackrel{\gamma_M}{\longmapsto} arm$ is an isomorphism
\begin{equation}
R \o_T (A \o_B M) \stackrel{\cong}{\longrightarrow} M
\end{equation}
of left $A$-modules.  
\end{lemma}
\begin{proof}
Note that the left $A$-module structure on $R \o_T A \o_B M$ is given by
$a \cdot \, r \o_T x \o_B m = r \o ax \o m$  ($a, x \in A$), which is well-defined since
since the left $T$-module structure on the induced module
$A \o_B M$ is given by $t \cdot a \o m = at^1 \o t^2m$.  
An inverse mapping  $M \to R \o_T (A \o_B M)$ is defined in terms of a separability element by
$m \mapsto 1_A \o e^1 \o e^2m$ for each $m \in M$, since on the one hand we have $e^1 e^2 m = m$
and on the other hand, given $a \in A, r \in R$, 
$$ 1 \o (e^1 \o e^2 arm) = 1 \o_T (e^1 \o e^2r) \cdot (a \o m) = r \o (a \o m). \qed
$$
\renewcommand{\qed}{}\end{proof}

The lemma recovers the isomorphism in \cite[Theorem 4.1]{KK} when $M = A$.  
Let $F$ denote the induction functor on $A$-modules, so that $F(M) = A \o_B M$
for an $A$-module ${}_AM$ and $F(f) = \id_A \o f$ for an arrow $f \in {}_A\M$. 
In the theorem below, $F$ is shown to be a fully faithful functor between ${}_A\M$
and its image $F({}_A\M)$ (perhaps more traditionally denoted by $F({}_A\M) = \Ind \Res {}_A\M$)
\cite[p.\ 95]{MAC}. In addition, the lemma together with projectivity of $R_T$ is clearly
seen from the proof to characterize separability.  

\begin{theorem}
\label{th-sep}
A ring extension $A \| B$ is separable if and only if the module $R_T$ is  projective
and the induction functor $F: {}_A\M \to F({}_A\M)$
 is an equivalence of categories with inverse functor given by $G(A \o_B M) = R \o_T (A \o_B M)$.   
\end{theorem}
\begin{proof}
($\Rightarrow$)  That the module $R_T$ (given by $r \cdot t = t^1 r t^2$) is cyclic projective
follows from \cite[4.1]{KK}, or more simply by noting that the right $T$-epi 
\begin{equation}
\eps_T: T_T \to R_T: \ \ \eps_T(t) = t^1 t^2 
\end{equation}
is split by $r \mapsto e^1 \o e^2 r$ using a separability element $e$.

By lemma, we see that there is a natural isomorphism $\gamma:\, GF \to I$ given
by $\gamma_M(r \o_T a \o_B m) := arm$ on an object $M$, a natural transformation
since for each arrow $f: N \to M$ in ${}_A\M$, the equation $\gamma_M \circ GF(f) = f \circ \gamma_N$
just follows from $f$ being $A$-linear.
Similarly, we see that there is a natural isomorphism $\phi: FG \to I$
given by $ \phi_M: A \o_B (R \o_T (A \o_B M)) \stackrel{\cong}{\longrightarrow} A \o_B M$,  $\phi_M(a \o_B r \o_T x \o_B m) = a \o_B xrm$ on an object $M$ in ${}_A\M$, where for each arrow $\id_A \o g: A\o_B N \to A\o_B M$ in $F({}_A\M)$
we have $(\id_A \o g) \phi_N = \phi_M (\id_A \o \id_R \o \id_A g)$ since $g$ is $A$-linear.
Thus, $F$ is an equivalence of the category ${}_A\M$ with its image subcategory.

($\Leftarrow$) Let $M$ be the natural left module ${}_AA$.  Applying the hypothesis on the ring extension $A \| B$,
we have $R_T$ is  projective and $R\o_T (A \o_B A) \cong A$ as left $A$-modules, and  
by naturality, as $A$-$A$-bimodules.  By \cite[4.1]{KK} this characterizes a separable
extension (via an application of the functor $- \o_T (A \o_B A)$ to $R_T \oplus * \cong  T_T$
and the isomorphism $GF(A) \cong A$).
\end{proof}

For example, if $A \| B$ is a separable algebra, the centralizer $R = A$ and $T = A^e = A^{\rm op} \o A$.  The module $A_{A^e}$ is
projective, a characterization of separability \cite{DI}. The
epimorphism $\eps_T: T \to R$ corresponds to the well-known
epimorphism $\mu: A^e \to A$ with interesting kernel,
the module of universal differentials $\Omega^1 = \{ x(dy)| x,y \in A \}$ where 
$dx = x \o 1 - 1 \o x$ (cf.\ \cite{BB}).  

The following corollary is the well-known relative projectivity characterization of separable
extension; e.g. \cite[10.8]{P}.  From the commutative diagram we see that tensoring
by $R_T$ is the ingredient that changes split epis to isomorphisms in a study
of separability.

\begin{cor}
A ring extension $A \| B$ is separable if and only if the mapping $\mu_M : A \o_B M \to M$ given
by $\mu_M(a \o m) = am$ is a split $A$-$C$-epimorphism for each $A$-$C$-bimodule $M$ where $C$ is a ring.  
\end{cor}
\begin{proof} 
Let 
$A \| B$ be any ring extension and  $M$ any $A$-$C$-bimodule. 
The homomorphism $\psi_M: A \o_B M \to R \o_T A \o_B M$,
$\psi_M(a \o m) = 1 \o a \o m$, which is induced by the epi $\eps_T: T_T \to R_T$, leads us to a commutative
diagram:

$$\begin{diagram}
A \o_B M && \rTo_{\psi_M} && R \o_T (A \o_B M) \\
&\SE_{\mu_M} & & \SW_{\gamma_M} & \\
& &  M & & 
\end{diagram}$$

\vspace{.5cm}

If $A \| B$ is separable, then $\gamma_M$ is an isomorphism of $A$-$C$-bimodules
and $\eps_T$ is a split epi, whence $\psi_M$ is split.  It follows that $\mu_M$ is a split
$A$-$C$-epimorphism as stated.  

Of course, if $\mu_A$ is split as an $A$-$A$-epi, then $A \| B$ is separable by definition.
\end{proof}

  Next we recall that a ring extension $A \| B$ is \textit{split}, if $B \to A$ has a left inverse 
as $B$-$B$-homomorphisms; thus $A \| B$ is a proper extension,
$A$ is a right and left $B$-generator and a left inverse $E$, called
a \textit{conditional expectation}, satisfies $E(1_A) = 1_B$. 
Recall that for a ring homomorphism $B \to A$, the coinduced module $(\Coind N)_A$ of a $B$-module $N_B$
is an $A$-module $\Hom (A_B, N_B)$, which is of course a right $S$-module via
composition with elements of $S = \End {}_BA_B$. 
The module action of $S$ however only commutes with the restricted action of $A$
to $B$.

\begin{lemma}
\label{lem-split}
Given a split extension $A \| B$ and a $B$-module $N_B$,
the mapping $f \o r \mapsto f(r)$ defines an isomorphism
\begin{equation}
\Hom (A_B, N_B) \o_S R \stackrel{\cong}{\longrightarrow} N_B
\end{equation}
of right $B$-modules.
\end{lemma}
\begin{proof}
The mapping is $S$-balanced since ${}_SR$ is given by evaluation,
$\alpha \cdot r = \alpha(r)$.  The mapping is right $B$-linear since
the $B$-module structure is given by $(f \o r) \, \cdot b = fb \o r$, which
is $S$-balanced as well, 
and $f(br) = f(rb) = f(r)b$ for $f \in \Hom (A_B, N_B), b \in B,
r \in R = A^B$.  

An inverse is given $n \mapsto nE(-) \o 1_A$ where $n \in N$ and $nE(-)$ denotes
the mapping in $\Hom (A_B, N_B)$ given by $x \mapsto nE(x)$.  This is indeed an
inverse since $nE(1) = n$ on the one hand and 
$$  f(r)E(-) \o 1_A = f \circ \rho_r \circ E \o_S 1_A  = f \o r, $$
on the other hand.  
\end{proof}

Now consider the functor corresponding to the restriction of coinduction,
$F: \M_B \to F(\M_B)$ given by $F(N_B) = \Hom (A_B, N_B)_B$. (The image
of $F$ is the category we might denote by $\Res  \Coind  \M_B$.)
For a split extension this functor is an equivalence as we see next.

\begin{theorem}
\label{th-split}
A ring extension $A \| B$ is split if and only if the module ${}_SR$ is  projective
and the induction functor $F: \M_B \to F(\M_B)$
 is an equivalence of categories with inverse functor given by $G(\Hom (A_B, N_B)) = \Hom (A_B, N_B) \o_S R$.   
\end{theorem}
\begin{proof}
($\Rightarrow$) That ${}_SR$ is cyclic projective follows from \cite[Lemma 1.1]{LK}, or by other means,
noting that the left $S$-homomorphism $\eps_S: \, S \to R$ given by $\alpha \mapsto \alpha(1)$
splits the monomomorphism ${}_SR \to {}_SS$ given by
$r \mapsto rE(-)$. 

By lemma $GF \cong I$ on objects in $\M_B$, and is natural w.r.t.\ arrows $f: N_B \to U_B$ in $\M_B$ 
since the evaluation mappings in the lemma clearly commute with $\Hom (A,f) \o \id_R$ (where
$\Hom (A, f) = \lambda_f$) and $f$.
Also $I \cong FG$ on objects $\Hom (A_B, N_B)_B$ in $F(\M_B)$, since $$\Hom (A_B, N_B) \stackrel{\cong}{\longrightarrow}
\Hom(A_B, \Hom (A_B, N_B) \o_S R), \ \ g \longmapsto (a \mapsto g(a)E(-) \o 1_A) $$
and is natural w.r.t.\ an arrow $f: N_B \to U_B$, since $f(g(a)E(-)) = f(g(a)) E(-)$.  

($\Leftarrow$)  Let $N$ be the natural module $B_B$, whence $\Hom (A_B, B_B) \o_S R \cong B$
as right $B$-modules, and also as left $B$-modules by naturality.  Apply the
functor $\Hom (A_B, B_B) \o_S -$ (from left $S$-modules to $B$-$B$-bimodules)   
to ${}_SR \oplus * \cong {}_SS$.  Whence
$$ {}_BB_B \oplus * \cong  {}_B\Hom (A_B, B_B)_B. $$
It follows that there is $f \in \Hom ({}_B\Hom(A_B, B_B)_B, {}_BB_B)$
and $$g \in \Hom ({}_BB_B, {}_B\Hom (A_B, B_B)_B) \cong \Hom ({}_BA_B, {}_BB_B)$$ via
$g \mapsto g(1_B)$, such that $ f \circ g = \id_B$.
If $F := g(1)$, this becomes $ f(F) = 1_B$.  Define $E(a) = f (\lambda_a \circ F)$,
an arrow in $\Hom ({}_BA_B, {}_BB_B)$ satisfying $E(1_A) = 1_B$.
\end{proof}

Note that the lemma and the proof of the theorem show that
\begin{cor}
A ring extension $A \| B$ is split if and only if the module ${}_SR$ is finite projective
and $\Hom (A_B, B_B) \o_S R \cong B$ as $B$-$B$-bimodules.  
\end{cor}

%%%%%%%%%%%%%%%%%%%%%%%%%%%%%%%%%%%%%%%%%%%%%%%%%%%%%%%%%%%%%%%%%%%%%%%%%%%%%%%%%%%%%%%%%%%%%%%%%%%%%%%%%%%%%%%%
\section{Left D2 and H-separable Extensions}

A ring extension $A \| B$ is said to be \textit{left D2} (or left depth two) if 
its tensor-square $A \o_B A$ is centrally projective w.r.t.\ $A$ as natural $B$-$A$-bimodules; i.e.,
\begin{equation}
\label{eq: ld2mod}
{}_BA \o_B A_A \oplus * \cong \oplus^n {}_BA_A
\end{equation}
Under some natural identifications of $\Hom ({}_BA_B, {}_BA \o_B A_A) \cong T$ and
$\Hom ({}_BA \o_B A_A, {}_BA_A) \cong S$ valid for any ring extension, we arrive at
the equivalent condition of left D2:  there are matched elements $\beta_i \in S$
and $t_i \in T$ (called \textit{left D2 quasibases}) such that the following equation in $A \o_B A$ holds, 
\begin{equation}
\label{eq: left d2 qb}
x \o y = \sum_{i=1}^n t_i^1 \o t_i^2 \beta_i(x) y
\end{equation}
for all $x, y \in A$.

By letting $y = 1$, this last equation is quite analogous to the equation $a\1 \o \tau(a\2)a\3  = a \o 1$ for 
$a$ in some Hopf algebra with antipode $\tau$. Let's see this idea in action in the
next lemma, which is going to be used in our last section on normality.

\begin{lemma}
\label{lemma-norm}
Let $A \| B$ be a left D2 extension and ${}_A M_A$ a bimodule.  Then $AM^B \subseteq M^B A$.
\end{lemma}
\begin{proof}
For each $m \in M^B = \{ m \in M \, :\, mb = bm, \ \forall \, b \in B \}$, define
a mapping $\eps_m: A \o_B A \to M$ by $\eps_m(x \o y) = xmy$, which is well-defined.  
Extending the notion of the  module $R_T$, one defines a right $T$-module structure on $M^B$ in eq.~(\ref{eq: embeetee})
by $m \cdot t = t^1 m t^2$.   
Apply $\eps_m$ to eq.~(\ref{eq: left d2 qb}) with $y = 1$ to obtain $xm = \sum_i (m \cdot t_i) \beta_i(x) \in M^B A$.
\end{proof}

Naturally we would get equality in the lemma if we also had a corresponding right D2 condition, something we will do
in the next section (and such extensions are said to be \textit{D2}).  
For the next lemma we recall that we have $T \o_R A \cong A \o_B A$ via a
mapping $\pi_A$ sending $t \o_R a
 \mapsto t^1 \o_B t^2a$ with inverse given by $x \o_B y \mapsto \sum_i t_i \o_R \beta_i(x)y$
(an application of eq.~(\ref{eq: left d2 qb}) \cite[2.1(4)]{KK}). 

\begin{lemma}
\label{lem-d2}
If $A \| B$ is left D2 and ${}_AM$ is a module, then $\gamma_M(r \o a \o m) = arm$ defines an isomorphism
\begin{equation}
R \o_T (A \o_B M) \stackrel{\cong}{\longrightarrow} M
\end{equation}
of right $A$-modules, which is natural w.r.t.\ $A$-module homomorphisms ${}_AM \to {}_AN$.
\end{lemma}
\begin{proof}
Define a mapping $\pi_M =  \pi_A \o_A \id_A $ via the canonical identification
$A \o_A M \cong M$ as a component of the top horizontal arrow below. Since $A \| B$
is left D2, it follows that this arrow is an isomorphism.   
 The left
vertical and bottom horizontal arrows are again canonical isomorphisms.

$$\begin{diagram}
R \o_T T \o_R M &&\rTo^{\cong}_{\id_R \o \pi_M} && R \o_T (A \o_B M)\\
\dTo^{\cong} && && \dTo_{\gamma_M}\\
R \o_R M && \rTo^{\cong} && M 
\end{diagram}$$ 

\vspace{.5cm}
  
The diagram is commutative by following an element $r \o t \o m$ in the upper lefthand corner
around from below to $t^1 r t^2 m$ and around from above to $r \o t^1 \o t^2m$, then to $t^1 r t^2m$,
the same element.  Thus, $\gamma_M$ is an isomorphism. It is clearly
 natural w.r.t.\ arrows in ${}_A\M$ as in the proof of Theorem~\ref{th-sep}. 
\end{proof} 

Let $F$ again denote the induction functor on $A$-modules, so that $F(M) = A \o_B M$
for an $A$-module ${}_AM$. 

\begin{theorem}
\label{th-ld2}
If a ring extension $A \| B$ is left D2, then the induction functor $F: {}_A\M \to F({}_A\M)$
 is an equivalence of categories with inverse functor given by $G(A \o_B M) = R \o_T ( A \o_B M)$.   
\end{theorem}
\begin{proof}
By lemma, $GF \cong I$ via the natural transformation $\gamma_M: GF(M) \stackrel{\cong}{\longrightarrow} M$.
Also $FG \cong I$ on $F({}_A\M)$ by applying $F$ to the natural transformation $\gamma$.  
\end{proof}

In contrast to Theorem~\ref{th-sep} for separable extensions, 
I am not aware of any characterization of left D2 extensions in terms of $F$ and $G$ with some condition
on $R_T$ alone. There are examples that suggest we avoid weakening or  strengthening  the
projectivity condition on $R_T$, while the theorem below suggests weakening a generator condition
on $R_T$.    
There are of course separable extensions that are not left  D2 such
as the complex group algebras corresponding to a non-normal subgroup pair $H < G$ in a finite
group $G$ \cite[3.2]{KK}.  There are also left D2 extensions that are not separable 
such as the group algebras corresponding to a proper normal subgroup
 $H \ract G$ of a $p$-group over a characteristic $p$ field;
or more mundanely, any finite projective algebra that is not
a separable algebra.
The H-separable extensions (named after Hirata) are both separable, left (and right) D2 \cite[section~3]{KS}
for which we have the following theorem.  Recall that $A \| B$ is H-separable if there are matched elements
(called an H-sepability system) 
$e_i \in (A \o_B A)^A$ and $r_i \in R$ ($i = 1,\ldots,n$) such that in $A \o_B A$ we have $1 \o 1 = \sum_i e_i r_i$. For example,
from an H-separability system we see $A \| B$ is left D2 with
$\beta_i = \rho_{r_i}$ and $t_i = e_i$.   

\begin{theorem} 
\label{th-hsep}
A ring extension $A \| B$ is H-separable if and only if $R_T$ is a progenerator 
and $F$, $G$ defined in Theorems~\ref{th-sep} and~\ref{th-ld2} are inverse equivalences.
\end{theorem}
\begin{proof}
($\Rightarrow$) This follows from Theorems~\ref{th-sep} and~\ref{th-ld2} and the observation
that $R_T$ is a generator \cite[4.2]{KK}, which we also see as follows.  Given an H-separability system, the
epi $\oplus^n R_T \to T_T$ given by $(x_1,\ldots,x_n) \mapsto \sum_{i=1}^n x_i e_i$ is
split by the right $T$-monomorphism $t \mapsto (t^1r_1 t^2,\ldots, t^1 r_n t^2)$.

($\Leftarrow$) From the generator condition $T_T \oplus * \cong \oplus^n R_T$ we arrive
at $${}_AA \o_B A_A \oplus * \cong \oplus^n {}_AA_A,$$
which is another characterization of H-separable extension $A \| B$ \cite[ch.~2]{K}, by tensoring with $- \o_T (A \o_B A)$
and using naturality together with $GF \cong I$.
\end{proof}

In this theorem, $F$ and $G$ are not Morita equivalences, although if we expand the category $F({}_A\M)$
to the category ${}_T\M$, the functor $R \o_T -$ is a Morita equivalence to the category
${}_{Z(A)}\M$, where $Z(A)$ is the center of $A$ \cite[4.3]{KK}.
For example, if $A \| B$ is an Azumaya algebra (a separable algebra
over its center $B$), then we recover the well-known  fact that
$A^e$ and $B$ are Morita equivalent via the progenerator $A_{A^e}$
\cite{DI,K}.  

We now turn to coinduction of modules in relation to left D2 extension. The lemma below extends an aspect of the characterization
of the endomorphism ring as a smash product with the bialgebroid
$S$ in \cite[3.8]{KS}. 

\begin{lemma}
If the extension $A \| B$ is left D2 and $M_A$ is a module, then there is a isomorphism
$\Hom (A_B, M_B) \cong M \o_R S$
of right $B$-modules and right $S$-modules, which is natural with respect to arrows in $\M_A$. 
\end{lemma}
\begin{proof}
The isomorphism is given by 
$$\chi_M: \ M \o_R S \to \Hom (A_B, M_B), \ \ \ \chi_M(m \o \alpha)(a) :=  m\alpha(a)). $$
This is $B$-linear w.r.t.\ the right $B$-module $M \o_R S$ given by $m \o_R \alpha\, \cdot b$
$ = mb \o_R \alpha$ since $\chi_M(mb \o \alpha)(a) = \chi_M(m \o \alpha)(ba)$. It is also
right $S$-linear since $\chi_M(m \o \alpha \circ \beta)(a) = \chi_M(m \o \alpha)(\beta(a))$
for $\alpha, \beta \in S$.  

The transformation $\phi$ is natural w.r.t.\ $g  :  M_A \to N_A$ since 
$g(m\alpha(-)) = g(m)\alpha(-)$ for each $\alpha \in \End {}_BA_B$.  

An inverse mapping is given by $$\Hom (A_B, M_B) \to M \o_R S: \ \ \ f \mapsto \sum_i f(t^1_i)t^2_i \o_R \beta_i, $$
inverse by short computations using $\alpha(t^1)t^2 \in R$ and $\mu_M (f \o \id_A)$ applied
to eq.~(\ref{eq: left d2 qb}). 
\end{proof}

Let $R: \M_A \to \M_B$ denote the usual restriction (pullback or forgetful) functor from
$A$-modules to $B$-modules for a ring homomorphism $B \to A$.  The next theorem is like Theorem~\ref{th-split}
for split extensions in that restricted
coinduction is shown to be a fully faithful functor
but restricted this time to the subcategory $R(\M_A)$ of $\M_B$.
Again let $F$ be the functor of restricted coinduction from $\M_B$ into $\M_B$ given
by $F(N) = \Hom (A_B, N_B)_B$ for any module $N_B$.  

\begin{theorem}
\label{th-left d2 coind}
If $A \| B$ is left D2, then $F$ is a category equivalence from $R(\M_A)$ into
$FR(\M_A)$
with inverse functor $G$ given on a module $M_A$ by $$G(\Hom (A_B, M_B)_B) = \Hom (A_B, M_B) \o_S R.$$ 
\end{theorem}
\begin{proof}
Note that applying the lemma and associativity of tensor product gives us
\begin{equation}
\Hom (A_B, M_B) \o_S R \cong M \o_R S \o_S R \cong M \o_R R \cong M
\end{equation}
a compositon of isomorphisms from left to right giving the mapping $\rho_M: \ f \o_S r \mapsto f(r)$
for $f\in \Hom (A_B, M_B), r \in R$. Since $\rho$ is natural w.r.t.\ to
morphism $R(f)$ for arrows $f$ in $\M_A$, it follows that  $GF \cong I$ on $R(\M_A)$.
From this it follows that also $FG \cong I$ on $FR(\M_A)$.  
\end{proof}

In contrast to Theorem~\ref{th-split} characterizing split extensions, 
which will also have $F$ as an equivalence on the subcategory $R(\M_A)$ of $\M_B$, 
I am not aware of any characterization of left D2 extensions in terms of $F$ and $G$ with some condition
on ${}_SR$ alone.  There are left  D2 extensions that are not split such
as the endomorphism algebra over the $2 \times 2$ upper triangular algebra \cite[5.5]{LK2}. 
 There are of course split extensions that are not left D2
such as the complex group algebras corresponding to a non-normal subgroup
 $H < G$ of a finite group \cite[3.2]{KK}, or certain 
exterior algebras over their centers.
If $A \| B$ is split and centrally projective (so ${}_BA_B \oplus * \cong \oplus^n {}_BB_B$,
thus it is left D2), 
we may characterize these
by placing the progenerator condition on ${}_SR$ \cite[1.1]{LK}
in a theorem formulated like Theorems~\ref{th-left d2 coind} and~\ref{th-split}. 

We end this section with an endomorphism ring characterization of left D2 for right f.g.\ projective 
extensions. Let $E = \End A_B$ denote the right endomorphism ring of a ring extension $A \| B$.  
This has $A$-$A$-bimodule structure given by $(xfy)(a) = x f(ya)$, equivalently
 $x \cdot f \cdot y = \lambda_x \circ f \circ \lambda_y$.  

\begin{prop}
Suppose $A \| B$ is a ring extension where $A_B$ is finite projective.  Then
$A \| B$ is left D2 if and only if 
\begin{equation}
\label{eq: endo}
{}_AE_B \oplus * \cong \oplus^n {}_AA_B.
\end{equation}  
\end{prop}
\begin{proof}
We note that for any ring extension, there is an isomorphism $\Phi$, 
\begin{equation}
\label{eq: h}
\End A_B \stackrel{\cong}{\longrightarrow} \Hom (A \o_B A_A, A_A),\ \ \ \Phi(f)(x \o y) = f(x)y,
\end{equation}
with inverse given by $F \mapsto F(- \o 1_A)$. Note that $\Phi$ is an $A$-$A$-isomorphism
since $\Phi(xfy)(a \o c) = x f(ya) c = x \Phi(f)(ya \o c)$. 

($\Rightarrow$) If we are given that the tensor-square is centrally projective w.r.t.\ ${}_BA_B$, 
\begin{equation}
\label{eq: lefty}
{}_BA \o_B A_A \oplus * \cong \oplus^n  {}_B A_A,
\end{equation}
 we apply the functor
$\Hom ( -, A_A)$ from $B$-$A$-bimodules into $A$-$B$-bimodules to this, 
obtaining eq.~(\ref{eq: endo}). This only uses eq.~(\ref{eq: h})
and does not use the hypothesis on $A_B$.

($\Leftarrow$)
Since $A_B$ is finite projective, it follows that $A \o_B A_A$ and ${}_AE$ are finite projective
by applying $- \o_B A$ and ${}_A\Hom (-, A_B)$ to $A_B \oplus * \cong \oplus^n B_B$.  
Hence the module $A \o_B A_A$ is  reflexive, so $\Hom ({}_AE, {}_AA) \cong A \o_B A$.  
Then we apply ${}_B \Hom ({}_A - , {}_A A)_A $ to eq.~(\ref{eq: endo})
and obtain eq.~(\ref{eq: lefty}) as a consequence.
\end{proof}

%%%%%%%%%%%%%%%%%%%%%%%%%%%%%%%%%%%%%%%%%%%%%%%%%%%%%%%%%%%%%%%%%%%%%%%%%%%%%%%%%%%%%%%%%%%%%%%%%%%%%%%%%%%%%%%%%%%%%%%
\section{Normality for Subrings}

A ring extension $A \| B$ is \textit{right D2}
if there are $n$ matched elements $u_j \in T$ and
$\gamma_j \in S$ such that the following equation
in the tensor-square $A\o_B A$ holds,
\begin{equation}
\label{eq: right d2 qb}
x \o y = \sum_{j=1}^n x \gamma_j(y) u_j^1 \o u_j^2
\end{equation}
for all $x, y \in A$.  

Right D2 extensions are dual to left 
D2 extensions in a certain category of ring extensions,
via the notion of opposite ring \cite{LK2}.  While there is no theoretical reason given yet 
for right D2 extensions to  be left D2, there has never
been observed a one-sided D2 ring extension that is not two-sided.
We say a ring extension is \textit{D2} if it is both left and right D2. 
For example, finite right Galois extensions with Hopf algebroid actions are automatically right D2,
but they can be shown to be left D2 as well via the antipode.  
Finding a Galois extension without antipode might give us a one-sided
D2 extension; however, there is Schauenburg's theorem that an antipode
is definable from a Galois isomorphism for a bialgebra action. 

In this section we study normality for subrings in relation to D2 subrings (or 
proper extensions via inclusion). The reason for this study is that normal subgroups of finite groups correspond
precisely to D2 complex  finite dimensional subgroup algebras \cite{KK}, which  led to Nikshych asking the interesting generalized question
of whether a D2  Hopf subalgebra is normal; we  answer this question affirmatively in case
 the Hopf subalgebra enjoys equality with its double centralizer in corollary~\ref{cor-tys}.   
We use the simplest possible definition
of normal subring in \cite{R} while dropping the requirement that
the ring and subring be semisimple artinian.
(The requirements of semisimplicity and the ideals be maximal are not needed in our study; however, the reader is urged to read 
\cite{R} before making any further use of this notion.)   First,
let $B$ be a subring of $A$; we say an (two-sided) ideal $I$ in $B$
is \textit{$A$-invariant} if $AI = IA$ as subsets of $A$.
An ideal $J$ in $A$ may be \textit{contracted} to an
ideal
 $J \cap B$ in $B$.
   
\begin{definition}
\begin{rm}
For the purposes of this section,
a subring $B \subseteq A$ is said to be \textit{normal} if all contracted ideals
are $A$-invariant; i.e., for every ideal $J$ in $A$, we have $A(J \cap B) = (J \cap B)A$.    
\end{rm}
\end{definition}

For example, the center of a ring is  a normal subring.
Normal subgroups correspond to normal subrings via group algebra
by the proposition on Hopf subalgebras below.

\begin{prop}
\label{prop-nc}
If $A \| B$ is a D2 extension, then the centralizer $R$ is a normal subring in  $A$.
\end{prop}

\begin{proof}
By lemma~\ref{lemma-norm} and its right D2 dual, we have the following equality of subsets
in an $A$-$A$-bimodule ${}_AM_A$,
\begin{equation}
\label{eq: a-inv}
AM^B = M^BA
\end{equation}
Let $J$ be an ideal in $A$.  Then $J^B = J \cap R$, so
$J$ contracted to $R$ is $A$-invariant by the equation.  
\end{proof}

For the same reasons as in eq.~(\ref{eq: a-inv}), each submodule $J \subseteq R_T$, given a D2 extension $A \| B$, satisfies
$A$-invariance:
\begin{equation}
AJ = JA
\end{equation}

 Recall that a Hopf subalgebra $K \subset H$ (with antipode $\tau$) is \textit{normal}
if $a\1 K \tau(a\2) \subseteq K$ and $\tau(a\1) K a\2 \subseteq K$
for all $a \in H$. The two terminologies do not conflict:

\begin{prop}
\label{prop-hop}
Suppose $K \subseteq H$ is a Hopf subalgebra.  Then
$K$ is normal subring of $H$ if and only if $K$ is a normal
Hopf subalgebra of $H$.
\end{prop}

\begin{proof}
($\Leftarrow$) Suppose $K$ is a normal Hopf subalgebra in $H$.  Let $J$
be an ideal in $H$ and $x \in J \cap K$.  Given $a \in H$,
we note that
\begin{equation}
ax = a\1 x \tau(a\2)a\3 \in (J \cap K)H
\end{equation}
by the normality condition.  Hence $H(J \cap K) \subseteq (J \cap K)H$, and similarly $(J \cap K)H \subseteq H(J \cap K)$.  

($\Rightarrow$) Here we use the well-known characterization of
normality for a Hopf subalgebra $K$ in $H$ by
$HK^+ = K^+ H$ \cite[ch.~3]{MON}where $K^+ = \ker \eps|_K$ and $H^+ = \ker \eps$ are the augmentation ideals
of the counit (and its restriction to $K$).  We note
that $K^+ = H^+ \cap K$, whence by hypothesis
$HK^+ = K^+ H$.  
\end{proof}

It is well-known in textbooks on group theory
that given a subgroup $H < G$, the centralizer  subgroup
$C_G(H)$ is normal in the normalizer subgroup $N_G(H)$;
consequently, if $H$ is normal in $G$, then $C_G(H)$ is
normal in $G$.  The converse though is not true by some
easy examples; e.g., within the dihedral group of transformations of a square. 
Here is a Hopf subalgebra analogue.

\begin{prop}
\label{prop-crop}
Suppose $K$ is a Hopf subalgebra of a finite dimensional
Hopf algebra $H$, in which the centralizer $C_H(K)$ is a Hopf subalgebra as well. If $K \subseteq H$ is a normal Hopf subalgebra, then
$C_H(K) \subseteq H$ is a normal Hopf subalgebra.
\end{prop}
\begin{proof}
Since $K \subseteq H$ is a normal Hopf subalgebra, it
is Hopf-Galois (with coaction by the quotient Hopf algebra
$H/ HK^+$) and therefore D2 \cite{KK}.
By proposition~\ref{prop-nc}, the centralizer $C_H(K)$ is a
normal subring in $H$, therefore a normal Hopf subalgebra by
proposition~\ref{prop-hop}. 
\end{proof}

%In the same spirit we supply a new proof to a special case of D2 Hopf subalgebras that are 
%normal: a finite Hopf-Galois $W$-extension $K \subseteq H$ is normal (where $W$ is a third Hopf algebra
%coacting on $H$ with coinvariants $K$).
%
% \begin{prop}\cite[3.1]{LK3}
%Suppose $K \subseteq H$ is a Hopf subalgebra, $H$ and $W$ are finite dimensional
%Hopf algebras such that $H$ is a $W$-extension of $K$ which is right Galois.
%Then $K$ is a normal Hopf subalgebra of $H$.
%\end{prop}
%\begin{proof}
%As $W$ has a right coaction on $H$, dually $W^*$ acts from the left on $H$ with
%invariant subalgebra $K$.  Montgomery \cite[8.3.1]{MON} implies that
%any left ideal $J$ that is stable under the action of $W*$ in $H$ satisfies $J = H(J \cap K)$; similarly
%any stable right ideal $I$ in $H$ satisfies $I = (I \cap K)H$. But $H^+$ is such a stable left
%and right ideal, since  

Notice that the proof shows the following.

\begin{cor}
\label{cor-tys}
Suppose $K$ is a Hopf subalgebra of a finite dimensional
Hopf algebra $H$, in which the centralizer $C_H(K)$ and the double centralizer $C_H(C_H(K))$ are Hopf subalgebras. If $K \subseteq H$ is a D2 Hopf subalgebra, then the double centralizer 
$C_H(C_H(K)) \subseteq H$ is a normal Hopf subalgebra.
\end{cor}

While $K \subseteq C_H(C_H(K))$ is always the case, this will be a proper subset
even for certain normal subgroups, such as the index $2$ subgroups in the quaternion $8$-element
group.  

%%%%%%%%%%%%%%%%%%%%%%%%%%%%%%%%%%%%%%%%%%%%%%%%%%%%%%%%%%%%%%%%%%%%%%%%%%%%%%%%%%%%%%%%%%%%%%%%%%%%%%%%%%%%%%%%%%

\section{Discussion}

In this section, we append some additional thoughts about the previous
four sections.  First,  we pose in terms of modules the problem of when a ring extension
$A \| B$ satisfies $R \o_T (A \o_B A) \cong A$ via the mapping $\gamma_A$ defined in section~2.  For
any ring $C$ and  modules $M_C, N_C$, this question becomes simply 
when is evaluation $\Hom (M_C, N_C) \o_E M_C \stackrel{\cong}{\longrightarrow} N_C$ an isomorphism where 
$E := \End M_C$?  We view this question in terms of Morita theory
and its generalizations, thereby extending Theorem~\ref{th-sep}.  

Secondly, two characterizations of the left D2 condition on a ring extension $A \| B$ are obtained via naturally isomorphic functors of
induction or coinduction from the category of $A$-modules into
the category of $B$-modules.  Third,  the centralizer
of a D2 extension is noted to be a pre-braided commutative ring.
Finally,  a depth two arrow in a bicategory of bimodules
is considered in comparison with the endomorphism ring extension
of a bimodule, and two candidates for a definition of D2 bimodule
are shown to be the same if the bimodule is finite projective.

\subsection{The Morita viewpoint on $\gamma_A$}: $R \o_T (A \o_B A) \stackrel{\cong}{\longrightarrow} A$.  Let $A \| B$ be a ring extension. An interesting question
is when $A \| B$ has isomorphic $\gamma_A$ (which recall is
given by $\gamma_A(r \o a \o c) = arc$ for $r \in R$, $a,c \in A$).  
Such an extension is  a generalization of separable extension,
right D2 extension and left D2 extension, which is clear from
lemmas~\ref{lem-sep} and~\ref{lem-d2} (as well as an adaptation of the
latter to
right D2 extensions). 

Thinking in terms of the canonical isomorphisms introduced in the
preliminaries of section~1, there are $R \cong {\rm Hom}_{A^e}(A \o_B A, A)$, $T \cong {\rm End}_{A^e} (A \o_B A)$, where $A^e = A^{\rm op} \o A$ and $\gamma_A$ is identical
with evaluation,
$$ \Hom_{A^e}(A \o_B A, A) \o_{{\rm End}_{A^e} (A \o_B A)}   (A \o_B A)
\longrightarrow A,$$
where $R_T$ is identical with the right ${\rm End}_{A^e}(A \o_BA)$-module ${\rm Hom}_{A^e}(A \o_B A, A)$ given by
right composition.
This leads to the extended question of when given a ring $C$,
there are modules $M_C$ and $N_C$ such that the evaluation mapping
\begin{equation}
\gamma_{M,N}: \ \ \Hom (M_C, N_C) \o_{\End M_C} M_C \stackrel{\cong}{\longrightarrow} N_C
\end{equation}
given by $\gamma_{M,N}(f \o m) = f(m)$, is an isomorphism? 

A solution to this when $N = C$ and $M_C$ is a progenerator
comes from Morita theory, since $\Hom (M_C, C_C) \o_E M_C \cong C$
as $C$-$C$-bimodules where $E := \End M_C$ (in addition to $M \o_C \Hom (M_C, C_C) \cong
\End M_C$ via rank one projection).

More generally, if given a $C$-module $X_C$, let $\sigma[X]$ denote
the full subcategory of right $C$-modules that are submodules of quotients 
of direct sums $X^{(I)}$ for an arbitrary indexing set $I$ \cite{W}.

\begin{prop}
Suppose $M_C$ is $\sigma[X]$-projective and $\sigma[X]$-generator
and $N_C \in \sigma[X]$ w.r.t.\ a third $C$-module $X$.
Then the evaluation mapping $\gamma_{M,N}$ above is a
right $C$-isomorphism.
\end{prop}   

The proof follows directly from the proof in \cite[p.\ 499]{Lam}.

If we pass to the Dress subcategory $\mathcal{D}[M]$ of $\sigma[M]$, whose objects
are direct summands of direct sums of copies of $M$, there is a 
solution given by the following.

\begin{theorem}
Suppose $M_C, N_C$ are two $C$-modules such that $N_C \in \mathcal{D}[M]$.  Then $\gamma_{M,N}$ is an isomorphism.  
\end{theorem}
\begin{proof}
Let $\pi_i: M_C \to N_C$  and $\iota_i: N_C \to M_C$ be
$2n$ mappings such that $\sum_{i=1}^n \pi_i \circ \iota_i = \id_N$.
Then an inverse to $\gamma_{M,N}(f \o m) = f(m)$ is given
by $n \mapsto \sum_i \pi_i \o \iota_i(n)$. 
\end{proof}

Lemma~\ref{lem-sep} is a corollary of this theorem, since
if $A \| B$ is separable, we have $N_C \oplus * \cong M_C$ via
the multiplication mapping,
where $C = A^e$, $N = A$ and $M = A \o_B A$.

As a last remark note that there is a version of the commutative
triangle in Figure~1 if ${}_DM_C$ is a bimodule w.r.t.\ rings $C$ and $D$.  This triangle consists of the evaluation mapping $\gamma_{M,N}$ and
the counit of adjunction $\eps$ to the standard Hom-Tensor 
adjunction between the functors $\Hom (M_C,-)$ and $-\o_D M$ between   $\mathcal{M}_C$ and $\mathcal{M}_D$, which
is another evaluation mapping. Between them is a canonical quotient mapping 
$\psi_{M}$ induced from left module structure $D \to \End M_C$. 
The commutative diagram is the
following.

 $$\begin{diagram}
\Hom (M_C, N_C) \o_D M&& \rTo_{\psi_{M,N}} && \Hom (M_C, N_C) \o_{\End M_C} M \\
&\SE_{\eps} & & \SW_{\gamma_{M,N}} & \\
& &  N & & 
\end{diagram}$$

\vspace{.5cm}

\subsection{Two Functorial Characterizations of Left D2 Extension.}
Given a ring extension $A \| B$, we consider two functors of induction
from the category of modules ${}_A\M$ into ${}_B\M$.  First we have
$I^A_B := \Res^A_B \Ind^A_B \Res^A_B$ which is defined as the notation
suggests, viz.\ $I^A_B({}_AM) = {}_BA \o_B M$.  Second,
we make use of constructed ring $T$ and its right module structure
over the centralizer $R$ given by 
$t \cdot r = t^1 \o t^2 r$ (recalled from the preliminaries of
section~1).  From this we obtain a second functor $I^T_R: {}_A\M \to
{}_B \M$ defined by $I^T_R({}_AM) = T \o_R M$ where the left $B$-module structure makes use of elements of $R$ and $B$ commuting
in $A$: $b \cdot (t \o m) = t \o bm$.  Left D2 extensions are characterized by these two functors being naturally isomorphic and having a special property.

\begin{theorem}
A ring extension $A \| B$ is left D2 if and only if $T_R$ if f.g.\
projective and the two functors $I^A_B$ and $I^T_R$ are naturally
isomorphic.
\end{theorem}
\begin{proof}
The theorem is a restatement of \cite[Prop.\ 2.2]{LK2} in terms of
functors.  If $A \| B$ is left D2, then $A \o_B M \cong T \o_R M$
via $a \o_B m \mapsto$ $\sum_i t_i \o_R \beta_i(a)m$ with inverse
$t \o_R m \mapsto t^1 \o_B t^2m$ in terms of a left D2 quasibase.
Either of these maps is clearly natural w.r.t.\ an arrow ${}_AM \to {}_AN$ in ${}_A\M$. That $T_R$ is finite projective follows from
applying eq.~(\ref{eq: left d2 qb}) to elements in $T = (A \o_B A)^B
\subseteq A \o_B A$.  

Conversely, given $T_R \oplus * \cong \oplus^n R_R$, apply the functor $ - \o_R A$ into ${}_B\M_A$, one obtains from naturality
and $I^A_B(A) \cong I^T_R(A)$ that 
${}_BA \o_B A_A \oplus * \cong \oplus^n {}_BA_A$, the left D2 condition.      
\end{proof}
 
Note that if $A \| B$ is left D2 then the two functors are also naturally isomorphic
as functors to the category of left bimodules over $T$ and $B$:
\begin{equation}
I^A_B \cong I^T_R: {}_A\M \to {}_{T \o B}\M
\end{equation}
where ${}_{T \o B}T \o_R M$ for a module ${}_AM$ is naturally given
by $(t \o b)\cdot (u \o m) = tu \o bm$,
and ${}_{T \o B}A \o_B M$ is given by
$(t \o b) \cdot (a \o_B m) = ba t^1 \o_B t^2m$. 

Now recall the endomorphism ring $S = \End {}_BA_B$
and introduce a third functor of coinduction
$CoI^S_R: {}_A\M \to {}_B\M$ defined by $CoI^S_R({}_AM) = {}_B\Hom ({}_RS, {}_RM)$  where one again makes use elements of $B$ and $R$ commuting
in $A$ in order to define the left $B$-module $\Hom ({}_RS, {}_RM)$.
In analogy with the previous results, the two functors below
are naturally isomorphic if $A \| B$ is left D2, also as
functors into ${}_{S \o B}\M$.  

\begin{theorem}
A ring extension $A \| B$ is left D2 if and only if the module
${}_RS$ is f.g.\ projective and $I^A_B$ is naturally isomorphic
to $CoI^S_R$.
\end{theorem}
\begin{proof}
This is a restatement of \cite[Prop. 2.2(ii)]{LK2} in functorial
terms.  If $A \| B$ is left D2, then $A \o_B M \cong \Hom ({}_RS, {}_RM)$ via \begin{equation}
\label{eq: chi}
\chi_M(a \o_B m) (\alpha) := \alpha(a)m
\end{equation}
for $\alpha \in S$, 
with inverse $F \mapsto \sum_i t_i^1 \o_R t_i^2 F(\beta_i)$.
Either mapping is natural w.r.t.\ an arrow $g: {}_AM \to {}_AN$
(where $\chi_M$ and $\chi_N$ fit into a commutative square with $\Hom (S, g)$ and
$A \o g$).  

Conversely, from ${}_RS \oplus * \cong \oplus^n {}_RR$,
apply ${}_B\Hom (- , {}_RA)_A$ to obtain 
via naturality and ${}_BI^A_B(A)_A \cong {}_BCoI^S_R(A)_A$
the left D2 condition, eq.~(\ref{eq: ld2mod}).  
\end{proof}

Note that $\chi_M$ is also a left $S$-homomorphism where
${}_SA \o_B M$ is naturally given by $\alpha \cdot a = \alpha(a)$
and ${}_S\Hom ({}_RS, {}_RM)$ is induced from the natural
module $S_S$.  This follows from $\chi_M(\beta(a) \o m) (\alpha)
= \chi_M(a \o m) (\alpha \circ \beta) $ for $\alpha, \beta \in S$.  

\subsection{The centralizer of a D2 extension is a pre-braided commutative ring.} In \cite{LK} a natural coaction $\delta: A \to A \o_R T$ on a (right) D2 extension $A \| B$
is introduced, and it is noted that $A \| B$ is a right $T$-Galois
extension if $A_B$ is faithfully flat (although this condition is
not needed for the discussion to follow).  The coaction $\delta$ is given 
in Sweedler notation by
\begin{equation}
a\0 \o a\1 := \sum_j \gamma_j(a) \o_R u_j
\end{equation}
 using the notation in eq.~(\ref{eq: right d2 qb}).  The coaction $\delta$ 
restricts to the centralizer $R$ of $A \| B$ to give the
mapping $s_R: R \to T$, where $s_R(r) = 1 \o_B r$,
 after the simplification $R \o_R T \cong T$:
$$ r\0 \o r\1 = \sum_j 1 \o \gamma_j(r)u_j = 1 \o_R 1 \o_B r. $$

Recall the module $R_T$ with action temporarily denoted by
$r \ract t = t^1 r t^2$.  Then the ring $R$ satisfies the following
\textit{pre-braided commutativity} condition (cf.\ \cite{BS}): ($\forall \, s,r \in R$)
\begin{equation}
sr = r\0 (s \ract r\1).
\end{equation}
Another equivalent way to look at this type of commutativity condition
is for a ring $R$ with module structures ${}_SR$ and $R_T$ which
do not combine to give a bimodule, but instead there
are elements $\gamma_j \in S$, $u_j \in T$ such that
\begin{equation}
sr = \sum_j (\gamma_j \cdot r)(s \ract u_j)
\end{equation}
The Miyashita-Ulbrich module $R_T$, the right $R$-bialgebroid $T$
and the right $T$-comodule $R$ formally satisfy the Yetter-Drinfeld
condition.    

\subsection{Depth Two Bimodules.} We consider possible definitions of a 
depth two bimodule from two viewpoints.  First, there is the viewpoint
of a depth two arrow in a bicategory with a recipe for right D2 in
\cite{S}.  There is a bicategory of rings where arrows or 1-cells are bimodules
and 2-cells are bimodule homomorphisms \cite[p.\ 283]{MAC}.  (Horizontal) composition of arrows is 
tensoring bimodules, which is only associative up to an isomorphism, while (vertical) composition of 2-cells is composition of bimodule homomorphisms.
If $\alpha: {}_RM_S \to {}_RN_S$, $\beta: {}_RN_S \to {}_RQ_S$, $\gamma: {}_SU_T \to {}_SV_T$, and $\delta: {}_SV_T \to {}_SW_T$ 
are bimodule homomorphisms, then the horizontal and vertical compositions satisfy a middle four exchange law
as an equality of homomorphisms between ${}_RM \o_S U_T \to {}_RQ \o_S W_T$
\begin{equation}
(\beta \o_S \delta) \circ (\alpha \o_S \gamma) = ( \beta \circ \alpha) \o_S (\delta \circ \gamma)
\end{equation}

To this we make some changes.  Replace an object ring $R$ by its category of modules $\M_R$,
then the arrows above become induction functors, to which we add coinduction functors $\Hom ({}_SM_R, -):
\M_R \to \M_S$.  The coinduction functor just given has left adjoint the induction functor $ - \o_S M_R: \M_S \to \M_R$.
The recipe for bimodule ${}_SM_R$ to be a right D2 arrow in \cite{S} is functorial with value on $M$ given by  
\begin{equation}
\label{eq: szlach}
\Hom (M_R, \End (M_R) \o_S M_R) \oplus * \cong \oplus^n \End (M_R) 
\end{equation}
as right $S$-modules.  For example, if $B \to A$ is a ring homomorphism and $M = {}_BA_A$ the induced bimodule,
then the induction functor is $\Ind^A_B$, the coinduction functor is $\Res^A_B$, and by naturality
one arrives at the right D2 condition,
\begin{equation}
\label{right d2 cond}
{}_AA \o_B A_B \oplus * \cong \oplus^n {}_AA_B.
\end{equation}

On the other hand, one might make use of the endomorphism ring  of a ring extension 
to define a bimodule ${}_SM_R$ to be right D2 if the left module structure mapping $S \to \End M_R$ is
right D2.  For example, given a ring homomorphism $B \to A$ and its induced bimodule ${}_BA_A$, this map simplifies
$B \to \End A_A \cong A$ to  $B \to A$. Now let $E$ denote $\End M_R$, then $S \to E$ is right D2 if
${}_EE \o_S E_S$ is centrally projective w.r.t.\ ${}_EE_S$.  What relation does the central projectivity
condition have to the condition in eq.~(\ref{eq: szlach})?  They are the same if $M_R$ is finite projective:

\begin{lemma}
Suppose ${}_SM_R$ is a bimodule which is f.g.\ projective as a right module.  
Then there is an $E$-$S$-bimodule isomorphism,
\begin{equation}
\End M_R \o_S \End M_R \stackrel{\cong}{\longrightarrow} \Hom (M_R, \End M_R \o_S M_R )
\end{equation}
given by $\phi(f \o g)(m) = f \o g(m)$ for every $f,g \in E$, $m \in M$.  
\end{lemma}

The proof of the lemma is left as an exercise.  
Note that for any ring extension $A \| B$ the canonical bimodule ${}_BA_A$ is right rank one free.    
These considerations lead
us then to define when a right finite projective bimodule is right D2.

\begin{definition}
A right f.g.\ projective bimodule ${}_SM_R$ is said to be right D2 if its endomorphism ring extension $S \to \End M_R$
is right D2.  
\end{definition}
  
The definition admits an extension of the definition to left D2 left finite projective bimodules by
opposite categorical considerations as in \cite{LK2}:  a left finite projective bimodule ${}_SM_R$ is left D2
if its left endomorphism ring extension $R \to (\End {}_SM)^{\rm op}$ is left D2.   Then a left D2 extension $A \| B$ yields
the left D2 bimodule ${}_AA_B$. However, left and right D2 bimodules
do not generalize properly D2 ring extension unless we restrict
ourselves to Frobenius extensions and Frobenius bimodules,
since a Frobenius extension is left D2 iff it is right D2.
Recall that a bimodule ${}_SM_R$ is Frobenius if ${}_SM$
and $M_R$ are finite projective modules, and the left and right
duals of $M$ are $R$-$S$-bimodule isomorphic: 
\begin{equation}
\Hom ({}_SM, {}_SS) \cong \Hom (M_R, R_R).
\end{equation}
The main property of a Frobenius bimodule is that its right
and left endomorphism ring extensions are Frobenius 
(Morita, cf.\ \cite[ch.\ 2]{K}).

\begin{definition}
A Frobenius bimodule ${}_SM_R$ is said to be D2 if
its endomorphism ring extensions $R \to (\End {}_SM)^{\rm op}$
and $S \to \End M_R$ are  D2.
\end{definition}

Note that a Frobenius extension $A \| B$ is D2 iff
its bimodule ${}_BA_A$ (or ${}_AA_B$) is Frobenius and D2,
which follows from the endomorphism ring theorem for D2
extensions \cite{LK2}.  
 For example, the Frobenius bimodule in \cite[2.7]{K} is D2
since its endomorphism ring is H-separable, therefore left and right D2.    
D2 bimodules would in principle provide a wealth of examples of  D2 extensions
via the endomorphism ring.

%%%%%%%%%%%%%%%%%%%%%%%%%%%%%%%%%%%%%%%%%%%%%%%%%%%%%%%%%%%%%%%%%%%%%%

\end{document}